\documentclass[12pt,a4paper,fleqn]{article}

\usepackage[centertags]{amsmath}
\usepackage{amsfonts}
\usepackage{amssymb}
\usepackage{amsthm}
\usepackage{amscd}
\usepackage{amsxtra}
\usepackage{newlfont}
\usepackage{titlesec}

\newcommand{\OK}{\mathfrak{O}_{K}}
\newcommand{\OF}{\mathfrak{O}_{F}}
\newcommand{\OKF}{\mathfrak{O}_{KF}}
\newcommand{\OKT}{\mathfrak{O}_{\widetilde K}}

\theoremstyle{plain}
\newtheorem{thm}{Theorem}
\newtheorem{cor}[thm]{Corollary}
\newtheorem{lem}[thm]{Lemma}
\newtheorem{prop}[thm]{Proposition}

\theoremstyle{definition}

\theoremstyle{remark}

\newtheorem*{remnn}{Remark}
\numberwithin{thm}{section}
\numberwithin{equation}{section}

\titleformat{\section}
{\normalfont\large\bfseries}{\S\thesection.}{1em}{}

\begin{document}
\pagestyle{myheadings}

\title{On the divisibility of class numbers of cubic number fields 
with discriminants in a prescribed rational quadratic class}
\author{Ivan Chipchakov\thanks{ The research of I. Chipchakov was supported
in part by Grant MM1106/2001 from the Bulgarian Foundation
for Scientific Research.}   and Kalin Kostadinov}
\date{}

\maketitle
\begin{abstract}
Let $n$ be an odd number and $F$ an imaginary quadratic field
with odd discriminant. We show that there exists infinitely 
many cubic fields $K$ such that the class number of $K$ 
is divisible by $n$ and the Galois closure of $K$ contains $F.$ 
\end{abstract}

\section{Introduction}
This paper is devoted to the study of the problem of 
whether there exists a cubic extension $K_{n,d}$ of 
the field $\mathbb{Q}$ of rational numbers such that 
the class number $\mbox{h}(K_{n,d})$ is divisible by 
an arbitrary fixed natural number $n$ and the discriminant 
$\mbox{disc}(K_{n,d})$ is lying in the set 
$\{ds^{2}:s\in\mathbb{N}\},$ 
for a given square-free integer $d.$
It is known that this problem has a positive solution 
for any $n$ in case $d=1$ \cite{uchida} 
and in case $d=-3$ \cite{nak3}.
Our main result shows that the answer to the considered 
question is affirmative whenever $n$ is odd and $d$ is 
negative and congruent to $1$ modulo $4.$
The field $K_{n,d}$ is associated with Uchida's solution 
$K_{n,1}$ by a congruence type relation modulo a suitably 
chosen set of prime numbers splitting in $K_{n,1}$ 
(see sections $\ref{sec:prsets}$ and $\ref{sec:mproof}$).

\section{Preliminaries on units and ramification in a 
cubic field generated by a root of Uchida's polynomial}
\label{sec:prelims}
\markright{\footnotesize\sc{CLASS NUMBERS OF CUBIC NUMBER FIELDS}}
Let $n,d$ and $a$ be integers such that $n$ is odd and positive, 
$d$ is negative, square-free and $d\equiv 1 \pmod{4}$, 
and $a$ is odd.
With these notations, the polynomial 
\[
 u(x)=x^{3}+m(x+1)^2, \quad
 \mbox{where } m=(3^{6}d^{n}a^{2n}+27)/4\, ,
\]
is irreducible over the field of rational numbers 
(cf. \cite[Lemma 1]{uchida}). 
Throughout this paper we denote by 
$K$ an extension of $\mathbb{Q},$
obtained by adjoining a root $\pi$ of the polynomial $u(x)$ 
and by $F$ the quadratic field 
$\mathbb{Q}(\sqrt{d})\;\;
(\pi, \sqrt{d}\in\overline{\mathbb{Q}},$ 
a fixed algebraic closure of $\mathbb{Q}$). 
It is easily verified that the discriminant 
$\mbox{disc}(u(x))$ of $u(x)$ 
is equal to $m^{2}(4m-27),$ 
i.e. the compositum $KF=K(\sqrt{d})$ 
is a root field of $u(x).$ 
As $\mbox{disc}(u(x))<0,\,K$ has only one embedding 
in the field of real numbers.
Hence $KF/\mathbb{Q}$ is an imaginary Galois extension of degree $6$,
and by Dirichlet$^{,}$s unit theorem, the multiplicative group 
$\OKF^{\ast}$ is of rank $2,\,\OKF$ being as usual 
the ring of algebraic integers in $KF.$
The subgroup of $\OKF^{\ast}$ singled out by the following lemma 
will play an essential role in our further considerations.
\begin{lem}\label{lem:eps}
Let $\sigma$ be an $F$-automorphism of the field $KF$ of order $3$.
Then the elements 
$\varepsilon =\pi+1$ and $\varepsilon^{\sigma} =\pi^{\sigma}+1,$
are units in $\OKF$ and generate a 
subgroup $E$ of finite index in $\OKF^{\ast}.$
\end{lem} 
\proof Clearly $\varepsilon$ and $\varepsilon^{\sigma}$ 
are roots of the polynomial 
$u(x-1)=x^{3}+(m-3)x^{2}+3x-1,$ 
i.e. they are elements of $\OKF^{\ast}.$ 
We prove that 
$[\OKF^{\ast}:E]\in\mathbb{N}$ 
by showing that $\varepsilon$ and $\varepsilon^{\sigma}$
are independent. Assuming that 
$\varepsilon^{s}=(\varepsilon^{\sigma})^{t},$ 
for some integers $s$ and $t,$ 
one obtains that 
$\varepsilon^{s}\in K\cap K^{\sigma}=\mathbb{Q}.$ 
As $\varepsilon^{s}\in\OKF^{\ast},$
this implies that $\varepsilon^{s}=\pm 1,$ 
which contradicts the fact that $1$ and $-1$ 
are all roots of unity in $K.$ \qed

The following proposition gives information 
about the discriminant $\mbox{disc}(K)$ 
and the ramification of primes not equal to $3$ 
in the ring $\OKF;$ 
the excluded case of $p=3$ 
is considered in the next section.

\begin{prop}\label{prop:lema1}
For a prime number $p$ not equal to $3$, 
the following is true:
\\
{\normalfont (i)} 
$p$ is totally ramified in $K$ 
if and only if 
it divides $m$ and $\mbox{disc}(K);$ 
this occurs 
if and only if 
the power of $p$ in the canonical 
primary decomposition of $m$ 
is not divisible by $3$;
\\
{\normalfont (ii)} 
$p.\OK = \mathfrak{p}_{1}^{2}.\mathfrak{p}_{2}$, 
where $\mathfrak{p}_{1}$ and $\mathfrak{p}_{2}$ 
are different prime ideals of $\OK$, 
if and only if
the power of $p$ in the canonical 
decomposition of the number $4m - 27$ is odd. 
\end{prop}

\markright{\footnotesize\sc{I. CHIPCHAKOV AND K. KOSTADINOV}}

\proof
(i)Consider the presentation 
$m = b.c ^{3}$ 
with $b, c$ integers and $b$ cube-free. 
The element $\rho=\pi/c$ of $K$ 
has a minimal polynomial 
\[
g(x)=c^{-3}.u(cx)=x^3+b(cx+1)^2,
\]
 with discriminant $b^{2}(4m-27).$ 
Note that it is enough to prove 
the following three implications:
\\
1. If $p$ divides both $m$ and $\mbox{disc}(K)$,
   then $p$ divides $b;$
\\
2. If $p$ divides $b,$ 
   then $p$ is totally ramified in $K;$
\\
3. If $p$ is totally ramified in $K,$ 
   then $p$ divides both $m$ and $\mbox{disc}(K).$
\\
For the first, note that 
$\mbox{gcd}(\mbox{disc}(K),m)$ 
is a divisor of 
$\mbox{gcd}(b^{2}(4m-27),m)$, 
and the last is a divisor of $27b^{2}.$ 
From $p\not=3$ it follows, 
that $p$ divides $b.$ \\
For the second, 
let $p$ be a prime divisor of $b$ 
not equal to $3.$ 
If $p^{2}\nmid b$, 
then $\rho$ is a root of an Eisensteinian polynomial 
relative to $p$ and in case $p^{2}|b$ 
this is true for the minimal polynomial of $\rho^{2}/p.$ 
Hence in both cases $p$ is totally ramified in $K$ 
(cf. \cite[Ch. I, Theorem 6.1]{cf}).\\
The proof of implication $3$ relies on the fact that the 
assumption on $p$ ensures that it divides 
$\mbox{disc}(K),$ which in its turn divides 
the product $m^{2}(4m-27).$ 
Note also that the norms
$N^{K}_{\mathbb{Q}}(\pi+3)\mbox{ and } N^{K}_{\mathbb{Q}}(\pi+m-6)$ 
are equal to $27-4m$ and $(m-8)(27-4m),$ respectively. 
It is therefore clear that if $p$ divides $4m-27,$ 
then the prime ideal $\mathfrak{p}$ of $\OKF$ lying above $p$ 
must contain the elements  $\pi+3,\;\pi+m-6$ and $4m-27.$
Since $4(\pi+3-(\pi+m-6))+(4m-27)=9,$ 
this leads to the conclusion that $9\in\mathfrak{p}$ 
in contradiction with the inequality $p\not=3.$ 
These observations prove that $p$ is a divisor of $m.$
\\
(ii) It is clear from Proposition \ref{prop:lema1}(i) 
that 
$p\OK=\mathfrak{p}_{1}^{2}.\mathfrak{p}_{2}\,
(\mathfrak{p}_{1}\not=\mathfrak{p}_{2})$
if and only if 
$p$ divides $\mbox{disc}(K)$ and $4m-27.$ 
In view of Dedekind's discriminant theorem 
\cite[Ch. III]{lang} 
and the inequality $p>3$ this occurs 
if and only if 
$p|\mbox{disc}(K)$ and $p^{2}\nmid\mbox{disc}(K).$ 
This combined with the fact that 
$\mbox{disc}(u(x))/\mbox{disc}(K)$ 
equals the square of the index $[\OK:\mathbb{Z}[\pi]]$ 
(as additive groups), proves Proposition \ref{prop:lema1}(ii).
\qed

\markright{\footnotesize\sc{CLASS NUMBERS OF CUBIC NUMBER FIELDS}}


\section{The canonical decomposition of the principal ideal 
$\frac{\pi-\pi^{\sigma}}{3\pi}\OKF$}
\label{sec:prepar}

The purpose of this section is to show that the element 
$\alpha=\frac{\pi-\pi^{\sigma}}{3\pi},$ 
where $\sigma$ is an $F$-automorphism of 
the field $KF$ of order $3,$ 
is an algebraic integer and 
the ideal $\alpha.\OKF$ is 
an $n-$th power of an ideal in $\OKF.$ 

\begin{lem}\label{prop:minpol}
The minimal polynomial of the element 
$\alpha=\frac{\pi-\pi^{\sigma}}{3\pi}$ 
over the field $F$ is 
\[
h(x)
=x^{3}+\frac{(2m-9+\sqrt{4m-27})}{6}x^2
-\frac{2m-9+3\sqrt{4m-27}}{18}x
+\frac{\sqrt{4m-27}}{27}.
\]
The coefficients of $h(x)$ are integers in $F$ 
and therefore $\alpha$ is in $\OKF.$
\end{lem}
\proof
It is convenient first to compute 
the minimal polynomial $h_{1}(x)$ over $F$ 
of the element $\pi^{\sigma}/\pi.$ 
The following equalities are true
\begin{eqnarray*}
Tr_{F}^{KF}(\pi^{\sigma}/\pi)+Tr_{F}^{KF}(\pi/\pi^{\sigma})=
\frac{\pi^{\sigma}}{\pi}+\frac{\pi^{\sigma^{2}}}{\pi^{\sigma}}
+\frac{\pi}{\pi^{\sigma^{2}}}+\frac{\pi}{\pi^{\sigma}}
+\frac{\pi^{\sigma}}{\pi^{\sigma^{2}}}+\frac{\pi^{\sigma^{2}}}{\pi }= 
&& \\
=(\pi+\pi^{\sigma}+\pi^{\sigma^{2}})
(\frac{1}{\pi}+\frac{1}{\pi^{\sigma}}+
\frac{1}{\pi^{\sigma^{2}}})-3=
-m.(-\frac{u'(0)}{u(0)})-3=2m-3,  
&& \\
Tr_{F}^{KF}(\pi^{\sigma}/\pi)-Tr_{F}^{KF}(\pi/\pi^{\sigma})=
\frac{\pi^{\sigma}-\pi^{\sigma^{2}}}{\pi}
+\frac{\pi^{\sigma^{2}}-\pi}{\pi^{\sigma}}
+\frac{\pi-\pi^{\sigma}}{\pi^{\sigma^{2}}}= 
&& \\
=\frac{(\pi-\pi^{\sigma})(\pi-\pi^{\sigma^{2}})
(\pi^{\sigma}-\pi^{\sigma^{2}})}{\pi.\pi^{\sigma}.\pi^{\sigma^{2}}}
= \sqrt{m^{2}(4m-27)}/(-m)= \sqrt{4m-27}.
\end{eqnarray*}

From the above it follows that
\begin{eqnarray*}
Tr_{F}^{KF}(\pi^{\sigma}/\pi)=(2m-3+\sqrt{4m-27})/2, \\
Tr_{F}^{KF}(\pi/\pi^{\sigma})=(2m-3-\sqrt{4m-27})/2.
\end{eqnarray*}

\markright{\footnotesize\sc{I. CHIPCHAKOV AND K. KOSTADINOV}}

The second elementary symmetric polynomial 
in the three variables 
$\frac{\pi^{\sigma}}{\pi},\frac{\pi^{\sigma^{2}}}{\pi^{\sigma}},
\frac{\pi }{\pi^{\sigma^{2}}}$ 
is equal to 
$Tr_{F}^{KF}(\pi/\pi^{\sigma})$, 
and the norm  of $\pi^{\sigma}/\pi$ is equal to $1.$
Therefore
\[
h_{1}(x)=x^3-\frac{2m-3+\sqrt{4m-27}}{2}x^2
+\frac{2m-3-\sqrt{4m-27}}{2}x-1.
\]
The minimal polynomial of 
$\alpha$ is computed from 
$h(x)=-\frac{1}{27}h_{1}(1-3x).$ 
\\
So at the end one gets
\[
h(x)=x^{3}+\frac{(2m-9+\sqrt{4m-27})}{6}x^2
-\frac{2m-9+3\sqrt{4m-27}}{18}x
+\frac{\sqrt{4m-27}}{27}.
\]
The free term of $h(x)$ is in $\OF,$ 
since  $3^{6}$ divides $4m-27$.
Since $27$ exactly divides $m,$ 
the numerators of the coefficients 
of the first and second degree terms of $h(x)$ 
are exactly divisible by $9$ 
and simple parity checks proves that 
these coefficients are also in $\OF.$
So $h(x)$ is in $\OF[x]$ and besides,
the trace of  $\alpha$ is a multiple of three, 
but the trace of $\alpha.\alpha^{\sigma}$ is not. 
\qed\\

The following two lemmas determine 
the decomposition in $\OKF$ of the prime ideals of $\OF$ 
containing the norm of $\alpha$ over $F.$

\begin{lem}\label{prop:33}
{\normalfont (i)} 
If $\,3\nmid da,$ 
then $\alpha$ is not contained in any ideal of $\OKF$ 
with norm equal to a power of $3.$ 
\\
{\normalfont (ii)} 
If $3|d,$ 
then the following decompositions are true:
\begin{eqnarray*}
3.\OF=P_{3}^{2},  
\\
P_{3}.\OKF=\mathfrak{P}_{3}^{'}\mathfrak{P}_{3}^{''}\mathfrak{P}_{3}^{'''}, 
\end{eqnarray*} 
where $P_{3}$ and 
$\mathfrak{P}_{3}^{'},\mathfrak{P}_{3}^{''},
\mathfrak{P}_{3}^{'''}$
are prime ideals respectively in $\OF$ and $\OKF.$ 
The element $\alpha$ belongs exactly to 
one of the ideals 
$\mathfrak{P}_{3}^{'},\mathfrak{P}_{3}^{''},\mathfrak{P}_{3}^{'''}.$
\end{lem}
\proof
(i) This follows immediately from the fact 
that in that case the norm of $\alpha,$ 
$N^{KF}_{\mathbb{Q}}(\alpha)=d^{n}a^{2n}$ 
is relatively prime to $3;$ \\
(ii) In this case the ideal $3\mathbb{Z}$ ramifies in $\OF$, 
since $3$ divides $\mbox{disc}(F)$. 
Since $3$ does not divide the trace of $\alpha.\alpha^{\sigma},$ 
there are at least  two different ideals of $\OKF$ 
with norm a power of $3$ and 
at least one of them does not contain $\alpha.$ 
This, combined with the fact that 
$KF/F$ is a cyclic cubic extension, 
implies that $P_{3}$ splits completely in $\OKF.$ 
Note finally that if $\alpha$ belongs to 
two of the ideals of $\OKF$ lying over $P_{3},$ 
then $\alpha.\alpha^{\sigma}\in P_{3}.\OKF,$ 
which implies that the trace of $\alpha.\alpha^{\sigma}$ 
is a multiple of  $3.$ 
The obtained contradiction proves the 
concluding assertion of Lemma \ref{prop:33}.
\qed\\
\markright{\footnotesize\sc{CLASS NUMBERS OF CUBIC NUMBER FIELDS}}



\begin{lem}\label{prop:lema2}
Let $P$ be a prime ideal of $\OF$,
 which contains $4m-27$ and does not contain $3.$
Then:\\
{\normalfont (i)} 
The ideal $P$ splits completely in $\OKF$;
\\
{\normalfont (ii)} 
The element $\pi-\pi^{\sigma}$ belongs exactly to 
one of the ideals of $\OKF$ which lie above $P$;
\end{lem}
\proof
(i) Denote by $p$ the prime number satisfying 
the equality $p\mathbb{Z}:=P\cap \mathbb{Q}.$ 
It is clear from the choice of $P$ that 
$p|4m-27$ and $p\not=3.$  
Therefore, $p$ does not divide $m,$ and 
by Proposition \ref{prop:lema1}, $
p$ is not totally ramified in $K.$ 
As $[KF:K]=2,$ this means that 
$3$ does not divide the ramification index of $p$ in $\OKF,$ 
which amounts to saying that $P$ does not ramify in $\OKF.$
Note also that $P$ is not inert in $\OKF.$ 
Assuming the opposite, one obtains from the multiplicative property 
of inertia degrees in towers of finite extensions 
that then $p$ must be inertial in $\OK.$ 
This, however, contradicts the existence of at least 
two different prime ideals of $\OK$ containing $p$ 
(established in the proof of Proposition \ref{prop:lema1}), 
and so proves our assertion. 
 As $KF/F$ is a cubic cyclic extension, 
these observations indicate that $P$ splits in $\OKF,$ as claimed.
 \\
(ii) It is easily verified that 
$u(x)$ has the following reduction modulo $4m-27:$ 
\begin{equation}\label{eq:red}
u(x)\equiv (x+3)^{2}(x+m-6)\pmod{4m-27}
\end{equation}
Let $P.\OKF=
\mathfrak{P}.\mathfrak{P}^{\sigma}.\mathfrak{P}^{\sigma^{2}}$.
It is clear from (\ref{eq:red}) that 
$(\pi+3)^{2}(\pi+m-6)$ belongs to 
$\mathfrak{P}$,$\mathfrak{P}^{\sigma}$  and$\mathfrak{P}^{\sigma^{2}}$. 
Since these ideals are prime, 
each contains at least one multiple of the above product.
Suppose for a moment that $(\pi+3)$ is contained in $\mathfrak{P}$,
$\mathfrak{P}^{\sigma}$ and $\mathfrak{P}^{\sigma^{2}}$, i.e. in their
intersection $P.\OKF$.
\markright{\footnotesize\sc{I. CHIPCHAKOV AND K. KOSTADINOV}}
Then the trace of $\pi+3$ over $\OF$ must lie in $P.\OF$. 
Since $\mbox{Tr}^{KF}_{F}(\pi+3)=-m+9,$ this means that 
$P.\OF$ contains the element $4m-27+4(9-m)=9.$ 
This contradicts the fact that the norm of $P.\OF$ 
is relatively prime to $3$ 
and so proves that $(\pi+3)$ is contained in at most two of 
$\mathfrak{P}$,$\mathfrak{P}^{\sigma}$ and $\mathfrak{P}^{\sigma^{2}}.$
Suppose that  $(\pi+m-6)$ is contained in two of the ideals, 
for example, $(\pi+m-6)\in\mathfrak{P}\cap\mathfrak{P}^{\sigma}$. 
Then 
$(\pi+m-6)(\pi^{\sigma}+m-6)
\in\mathfrak{P}.\mathfrak{P}^{\sigma}.\mathfrak{P}^{\sigma^{2}},$ 
and therefore 
$\mbox{Tr}^{KF}_{F}((\pi+m-6)(\pi^{\sigma}+m-6))\in P.$
On the other hand, direct calculations show that 
$\mbox{Tr}^{KF}_{F} ((\pi+m-6)(\pi^{\sigma}+m-6))=
4^{-2}\{(4m-61).(4m-27)+81\}.$ 
The obtained result implies that $81\in P$, 
which again contradicts the fact that 
the norm of $P$ is relatively prime to $3.$
Thus one concludes that $\pi+m-6$ lies in 
exactly one of the ideals 
$\mathfrak{P},\mathfrak{P}^{\sigma},\mathfrak{P}^{\sigma^{2}}.$
Assuming as we can 
that $\pi+3$ is contained in 
$\mathfrak{P}$ and $\mathfrak{P}^{\sigma},$ 
we obtain consecutively that
$\pi+m-6\in\mathfrak{P}^{\sigma^{2}},
\pi-\pi^{\sigma}\in\mathfrak{P}^{\sigma}$ 
and 
$\pi-\pi^{\sigma}\not\in\mathfrak{P}\cup\mathfrak{P}^{\sigma^{2}},$ 
which completes our proof. \qed 

\smallskip

The main result of this section can be stated as follows:
\begin{prop}\label{prop:ntapol}
 The principal ideal 
$(\alpha)=\frac{\pi-\pi^{\sigma}}{3\pi}.\OKF$ 
is an $n$-th power in the semigroup of ideals of $\OKF,$ 
i.e. there exist an ideal $\mathfrak{A}$ such that 
$(\alpha)=\mathfrak{A}^{n}.$
\end{prop}
\proof    
Since the norm $\mbox{N}^{KF}_{F}(\alpha)=\sqrt{4m-27}/27$ 
is an $n$-th power in $F^{\ast},$ 
this can be deduced from 
Lemmas \ref{prop:33} and \ref{prop:lema2}.
\qed\\

\section
{Sets of prime numbers associated with Uchida's cyclic cubic fields}
\label{sec:prsets}

In this section we present the main features of Uchida's 
solution of the class number divisibility problem 
for cyclic cubic fields and prepare the basis 
for constructing similar solutions for the case 
pointed out in the introduction. 
Our main result in this direction is the following theorem:


\begin{thm}\label{thm:petset} 
Let $\widetilde K$ be 
an extension of $\mathbb{Q}$ in $\overline{\mathbb{Q}}$ 
obtained by adjoining a root $\tilde{\pi} $ 
of the polynomial 
$\tilde{u}(X) = X^{3}+\tilde{m}X^{2}+2\tilde{m}X+\tilde{m},$ 
where $\tilde{m}=(\tilde{a}^{2^{s}n}+27)/4$, 
for some positive  integers 
$\tilde{a}$, $s$ and $n$. 
Assume also that 
$d$ is a square-free integer not equal to $1$, 
$\gcd(n, 6) = \gcd(\tilde{a}, 2) = 1$, 
$\tilde{\sigma}$ is an automorphism of $\widetilde{K}$ of order $3$, 
the ideal $\tilde{\alpha}.\OKT$, where 
$\tilde\alpha = (\tilde\pi - \tilde\pi^{\tilde\sigma })/\tilde\pi $, 
is not an $l$-th power of a principal ideal of $\OKT$, 
for any $l$ in the set 
$\mathcal{P}_{2n}$ of prime divisors of $2n,$ 
and the subgroup of $\OKT^{\ast}$ generated by the elements 
$\tilde\varepsilon=\tilde \pi + 1$ and 
$\tilde\varepsilon^{\tilde\sigma}=\tilde\pi^{\tilde \sigma } + 1$ 
is of index relatively prime to $2n$.
\\
Then there exist infinitely many pairs 
$(q_{1}(l; i, j), \, q_{2}(l; i, j))$ 
of prime numbers satisfying the following conditions, 
for each $l\in\mathcal{P}_{2n}$ 
and each pair $(i,j)$ of nonnegative integers less than $l:$
\\
{\normalfont (i)} 
$q_{1}$ and $q_{2}$ split in $\widetilde K;$
\\
{\normalfont (ii)} 
$q_{1} \neq q_{2}$ and 
$d$ is a $2^{s}$-th power residue modulo $q_{1}q_{2};$
\\
{\normalfont (iii)} 
$3$ is a $2^{s}n$-th power residue modulo $q_{1}q_{2};$
\\
{\normalfont (iv)} 
$\tilde\alpha_{ij}=
\tilde\alpha\tilde\varepsilon^{i} \tilde\varepsilon^{\tilde\sigma j}$ 
is not an $l$-th power residue modulo $q_{1}$ 
(i.e. modulo some prime ideal of $\OKT$ containing $q_{1}$);
\\
{\normalfont (v)} 
$\tilde \varepsilon^{i}\tilde\varepsilon^{\tilde\sigma j}$ 
is not an $l$-th power residue modulo $q_{2};$
\\
{\normalfont (vi)} 
$q_{1}$ and $q_{2}$ do not divide $6\mbox{disc}(\tilde{u}(X)).$
\end{thm} 
\markright{\footnotesize\sc{CLASS NUMBERS OF CUBIC NUMBER FIELDS}} 
Let us note that the existence of a 
cyclic cubic field $\widetilde{K}$ 
satisfying the conditions of Theorem \ref{thm:petset}, 
for any pair of integers $(n,s),$ 
has been established by Uchida \cite{uchida}.
As a matter of fact, he has shown there that 
the fulfillment of these conditions guarantees that 
the class number of $\widetilde{K}$ is divisible by $2^{s-1}n.$
\\
To prove the theorem, we need the following lemma:
\begin{lem}\label{lem:ass}
 With assumptions and notations being 
as in Theorem {\normalfont \ref{thm:petset},} 
fix a triple $(l;i,j)$ of admissible integers, 
and put  
\[
K_{1}=\widetilde{K}(\delta, \lambda, \zeta),\quad 
K_{2}=\widetilde{K}(\mu,\mu_{\sigma},\zeta_{0}),\quad 
K_{3}=\widetilde{K}(\xi,\xi_{\sigma},\zeta_{0})
\] 
where 
$\delta, \lambda, \zeta\,\mu,\mu_{\sigma},
\zeta_{0},\xi,\xi_{\sigma}$ 
lie in $\overline{\mathbb{Q}},\,
\delta^{2^{s}}=d, 
\lambda^{2^{s}n}=3,
\mu^{l}=\tilde\alpha_{ij},
\mu_{\sigma}^{l}=\tilde\alpha_{ij}^{\tilde\sigma},
\xi^{l}=\tilde\varepsilon,
\xi^{l}_{\sigma}=\tilde\varepsilon^{\tilde\sigma},
\zeta$ and $\zeta_{0}$ 
are primitive roots of unity of 
degrees $2^{s}n$ and $l$ respectively. 
Then $K_{2}$ and $K_{3}$ are not subfields of $K_{1}.$
\end{lem}
\markright{\footnotesize\sc{I. CHIPCHAKOV AND K. KOSTADINOV}}
\proof
It is clearly sufficient to show that 
$K_{1}$ does not include the fields 
$L_{2}=K_{2}(\zeta)$ and $L_{3}=K_{3}(\zeta).$ 
Note first that 
$K_{2}$ and $K_{3}$ are normal extension of $\mathbb{Q}$ 
with non-abelian Galois groups. 
Indeed, the fulfillment of the conditions of Theorem \ref{thm:petset} 
ensures the equality 
$[\widetilde{K}(\mu):\widetilde{K}]=
[\widetilde{K}(\xi):\widetilde{K}]=l,$ 
whereas the degree $[\widetilde{K}(\zeta_{0}):\widetilde{K}]$ 
is even and divides $l-1$ 
by the elementary properties of cyclotomic extensions. 
This implies that 
$\zeta_{0}\not\in\widetilde{K}(\mu)$  and 
$\zeta_{0}\not\in\widetilde{K}(\xi),$ 
which leads to the conclusion that 
$\widetilde{K}(\mu)/\widetilde{K}$ and 
$\widetilde{K}(\xi)/\widetilde{K}$ are 
not normal extensions, 
and thereby reduces our assertion to an 
immediate consequence of Galois theory. 
Since the field 
$L_{1}=\widetilde{K}(\zeta)$ 
is abelian over $\mathbb{Q},$ 
the obtained result shows that 
$K_{2}$ and $K_{3}$ are not included in $L_{1}.$ 
Applying now Kummer theory to the extensions 
$L_{2}/L_{1}$ and $K_{1}/L_{1},$ 
one concludes that if $L_{2}\subseteq K_{1},$ 
then there exist integers $n_{1}$ and $n_{2}$ 
such that 
$\gcd(n_{1},n_{2},l)=1,\,
\tilde\alpha_{ij}^{n_1}\tilde\alpha_{ij}^{\sigma n_{2}}
\not\in L^{\ast l}_{1}$ and
$(\tilde\alpha_{ij}^{\sigma}\tilde\alpha_{ij}^{-1})^{n_{1}}
(\tilde\alpha_{ij}^{\sigma^{2}}(\tilde\alpha_{ij}^{\sigma})^{-1})^{n_{2}}
\in L^{\ast l}_{1}.$ 
Observing, however, that 
$N^{\widetilde{K}}_{\mathbb{Q}}(\tilde\alpha_{ij}^{n_1}
\tilde\alpha_{ij}^{\sigma n_{2}})$ 
lies in $\mathbb{Q}^{\ast l},$ 
one verifies directly that 
$(\tilde\alpha_{ij}^{n_1}\tilde\alpha_{ij}^{\sigma n_{2}})^{3}
\in L_{1}^{\ast l}.$ 
Since $\gcd(l,3)=1,$ this yields 
$\tilde\alpha_{ij}^{n_1}\tilde\alpha_{ij}^{\sigma n_{2}}
\in L_{1}^{\ast l},$  
a contradiction, proving that 
$L_{2}$  is not a subfield of $K_{1}.$ 
Replacing $\tilde\alpha_{ij}$ by $\tilde\varepsilon_{ij}$ 
and arguing in the same way, one concludes that 
$L_{3}$ is not a subfield of $K_{1}$ either.
\vspace{0.2in}\qed
\\
{\it\hspace{-0.34in} Proof} of Theorem \ref{thm:petset}:
Suppose that the fields $K_{1},K_{2}, K_{3}$ 
are as in Lemma \ref{lem:ass} and for each $i$ 
denote by $\mbox{Spl}(K_{i})$ 
the set of rational prime numbers 
splitting completely in $K_{i}.$ 
It is clear from the definition of $K_{1}$ 
and the inclusion $\widetilde{K}\subset K_{1}$ 
that the elements of $\mbox{Spl}(K_{1})$ 
satisfy conditions (i),(ii) and (iii). 
It follows from Kummer's theorem that the elements of 
$\mbox{Spl}(K_{1})\setminus\mbox{Spl}(K_{2})$ 
satisfy condition (iv) 
and the elements of 
$\mbox{Spl}(K_{1})\setminus\mbox{Spl}(K_{3})$ 
satisfy condition (v) of Theorem \ref{thm:petset} 
(in either case with at most finitely many exceptions). 
Note finally that $K_{1},K_{2}, K_{3}$ are 
normal extensions of $\mathbb{Q}$, 
so Bauer's theorem (cf. \cite{hasse}) 
and Lemma \ref{lem:ass} imply that 
$\mbox{Spl}(K_{1})\setminus\mbox{Spl}(K_{2})$ 
and $\mbox{Spl}(K_{1})\setminus\mbox{Spl}(K_{3})$ 
are infinite sets. 
These observations prove Theorem \ref{thm:petset}. \qed
\begin{cor}\label{cor:qlij}
 With assumptions and notations being as in 
Theorem {\normalfont\ref{thm:petset},} let 
$\Omega_{\nu} = 
\{q(l; i, j): l \in\mathcal{P}_{2n}, 
i, j \in \{0,1..., l - 1\}\},\,(\nu=1,2)$ 
be sets of prime numbers with the properties 
required by Theorem {\normalfont \ref{thm:petset}.}
Then there exists an integer number $a$ 
satisfying the system of congruences
\[ 
3^{6}d^{n}a^{2^{s}n}\equiv\tilde{a}^{2^{s}n} 
\pmod{q_{1}(l;i,j)q_{2}(l; i, j)}: 
\quad q_{\nu}(l; i, j)\in\Omega_{\nu}. 
\]
\end{cor}
\proof  In view of the Chinese remainder theorem, 
it is sufficient to prove the 
solvability of the congruence 
$3^{6}d ^{n}x ^{2 ^{s}n} \equiv \tilde a^{2^{s}n} \pmod{q},$ 
for an arbitrary element $q\in\Omega_{1}\cup\Omega_{2}.$ 
The  choice of $\Omega_{1}\cup\Omega_{2}$  
ensures the existence of integers $z_{1}$ and $z_{2}$ 
satisfying the congruences 
$z_{1}^{2 ^{s}n} \equiv 3 \pmod{q}$ and
$z_{2}^{2 ^{s}} \equiv d \pmod{q}.$ 
It is therefore clear that every solution of the 
linear congruence
$z_{1}^{6}z_{2}x\equiv\tilde{a} \pmod{q}$ 
is a solution of the required type. 
\qed


\section{Proof of the main result}
\label{sec:mproof}

Let $n$ be a positive integer and 
$d<0$ and $d\equiv1\pmod{4}.$ 
In this section we find an infinite sequence of 
cubic fields with discriminants lying in 
$d\mathbb{Q}^{\ast2}$ and with 
class numbers divisible by $n.$ 
This is realized by a construction preserving 
the main features of Uchida's cyclic cubic fields.


\begin{thm}\label{thm:mproof}
Assume that 
$\widetilde K$, $\tilde \pi $, $\tilde m$, $\tilde a$, 
$n$, $s$ and $\tilde \sigma $ 
satisfy the conditions of 
Theorem {\normalfont \ref{thm:petset},} 
$\mathcal{P}_{2n}$ is the set of prime divisors of $2n$, 
$\Omega_{1}$ and $\Omega_{2}$ are sets of 
prime numbers determined in
accordance with Corollary {\normalfont \ref{cor:qlij},} 
$d$ is a square-free negative integer 
congruent to $1$ modulo $4$, 
and $K$ is an extension of $\mathbb{Q}$ 
in $\overline{\mathbb{Q}}$ 
obtained by adjoining a root $\pi $ of 
the polynomial 
$f(X) = X ^{3} + mX ^{2} + 2mX + m$, 
where $m = (3 ^{6}d ^{n}a ^{2 ^{s}n} + 27)/4$, 
for some integer number $a$ 
satisfying the system of congruences 
described in Corollary {\normalfont \ref{cor:qlij}.} 
Put $\varepsilon = \pi + 1$, $F = \mathbb{Q}(\sqrt{d})$ 
($\sqrt{d} \in \overline{\mathbb{Q}}$), 
and 
$\alpha = (\pi - \pi ^{\sigma })/(3\pi ^{\sigma })$, 
where $\sigma $ is a fixed automorphism of 
the compositum $KF$ of order $3$. 
Then the following assertions are true:
\\
{\normalfont (i)} 
The subgroup of $\OKF^{\ast}$ generated by 
$\varepsilon $ and $\varepsilon ^{\sigma }$ 
is of (finite) index 
relatively prime to $n$ and not divisible by $4.$
\\
{\normalfont (ii)} 
The ideals $\alpha.\OKF$ 
and 
$\alpha^{\sigma}\alpha^{\sigma^{2}}\OK$ 
are perfect $n$-th powers of ideals of $\OKF$ and $\OK$, 
respectively; 
\\
{\normalfont (iii)} 
The class groups of $KF$ and of $K$ 
contain elements of order $n;$
\\
{\normalfont (iv)} 
The class numbers $Cl(KF)$ and $Cl(K)$ are 
divisible by $3^{t}n$, 
in case at least $6 + t$ prime numbers are 
totally ramified in $K$.
\end{thm}


\begin{remnn}
  It should be noted that the 
system of congruences described in corollary \ref{cor:qlij} 
has a set $\{a_{k}:k\in\mathbb{N}\}$ of integer solutions, 
such that each $a_{k}$ gives rise to a 
cubic field satisfying the conditions of Theorem \ref{thm:mproof} 
for which the numbers of totally ramified 
primes is at least equal to $k.$ 
This is implied by the Chinese remainder theorem 
and the existence of infinitely many primes $p$ 
for which one can find an integer $b_{p}$ 
so that $\gcd(3^{6}d^{n}b_{p}^{2n}-27,p^{2})=p.$
\end{remnn}
{\it\hspace{-0.37in} Proof} of Theorem \ref{thm:mproof}: 
Suppose first that $q$ is an arbitrary element of 
$\Omega_{1}\cup\Omega_{2}.$ 
Our assumptions show that $u(x)$ and $\widetilde{u}(x)$ 
have one and the same reduction modulo $q,$ 
so $\mbox{disc}(u(x))\equiv\mbox{disc}(\widetilde{u}(x))\pmod{q},$ 
which ensures that $\gcd(\mbox{disc}(u(x),q)=1.$ 
This enables one to deduce the following statements from 
Kummer's theorem and Theorem \ref{thm:petset}(i):
\par (5.1)(i) The rational prime $q$ splits in $KF;$
\par \hspace{0.37in}(ii) $u(x)\equiv\widetilde{u}(x)\equiv 
(x-c_{1})(x-c_{2})(x-c_{3})\pmod{q\mathbb{Z}[x]},$ 
for some $c_{1},c_{2},c_{3}\in\mathbb{Z}.$ 
\\
Since $\mbox{Gal}(M/\mathbb{Q})\,:\,M\in\{\widetilde{K},KF\}$ 
acts transitively on the set of prime ideals of 
$\mathfrak{O}_{M}$ containing $q,$ and 
$(c_{3}-c_{1})(c_{1}-c_{2})(c_{2}-c_{3})$ 
is not divisible by $q,$ 
statement (5.1) can be supplemented as follows:
\par (5.2)(i) For each bijection $\nu$ of the set 
$\{1,2,3\},$ there exists a prime ideal 
$\mathfrak{Q}_{\nu}$ of $\OKF$ 
containing all four elements $
q,\pi-c_{\nu(1)},\pi^{\sigma}-c_{\nu(2)},\pi^{\sigma^{2}}-c_{\nu(3)};$
\par \hspace{0.35in}(ii) The ideal 
$\mathfrak{q}_{i}:=
q.\mathfrak{O}_{\widetilde{K}}+
(\widetilde{\pi}-c_{i}).\mathfrak{O}_{\widetilde{K}}$ is prime, 
for each index $i.$ 
Conversely, every prime ideal of 
$\mathfrak{O}_{\widetilde{K}}$ containing $q$ 
equals $\mathfrak{q}_{i},$ for some $i.$ 
\\
It is clear from (5.2)(ii) that if 
$\mathfrak{\widetilde{q}}$ is a 
prime ideal of $\mathfrak{O}_{\widetilde{K}},$ 
such that $q\in\mathfrak{\widetilde{q}},$ 
then $\mathfrak{\widetilde{q}}$ contains 
the elements 
$\widetilde{\pi}-c_{\nu(1)},
\widetilde{\pi}^{\widetilde{\sigma}}-c_{\nu(2)},
\widetilde{\pi}^{\widetilde{\sigma}^{2}}-c_{\nu(3)},$ 
for some bijection $\nu$ of the set $\{1,2,3\}.$ 
As $q$ is an arbitrary element of $\Omega_{1}\cup\Omega_{2},$ 
this observation, combined with (5.2)(i), 
leads to the following conclusions, 
for each admissible triple $(l,i,j):$
\par (5.3)(i) $\varepsilon_{ij}
=\varepsilon^{i}\varepsilon^{\sigma j}$ 
is not an $l$-th power residue modulo 
any $q_{2}(l,i,j)\in\Omega_{2};$
\par\hspace{0.35in}(ii) $\alpha_{ij}=
\alpha.\varepsilon_{ij}$ is not an $l$-th power modulo 
any $q_{1}(l,i,j)\in\Omega_{1}.$
\\
Note that Theorem \ref{thm:mproof}(i) 
is equivalent to the statement that every unit $\eta\in\OKF^{\ast}$ 
of norm $1$ over $F$ can be presented in the form 
$\eta=\varepsilon^{i}\varepsilon^{\sigma j}\eta_{1}^{2n},$ 
for some $\eta_{1}\in\OKF^{\ast}$ 
and suitably chosen integers $i,j.$ 
Therefore, its validity can be deduced from (5.3)(ii) 
and Dirichlet's unit theorem. 
Similarly, it follows from (5.3)(ii) and 
Theorem \ref{thm:mproof}(i) that 
$\alpha.\OKF$ is not an $l$-th power of a principal ideal, 
for any $l\in\mathcal{P}_{2n},\,l\not=2.$ 
\\
Observing that 
$\pi^{\sigma}/\pi^{\sigma^{2}}\in\OKF^{\ast},$ 
one sees that  $N^{KF}_{K}(\alpha^{\sigma}).\OKF=(\alpha^{\sigma})^{2}\OKF,$ 
which  reduces Theorem (ii) to a 
consequence of Proposition \ref{prop:ntapol}.
Note also that the obtained results prove Theorem \ref{thm:mproof}(iii). 
Applying finally the main result of  \cite{rz}, 
one completes the proof of Theorem \ref{thm:mproof}(iv).
\qed\smallskip
\markright{\footnotesize\sc{CLASS NUMBERS OF CUBIC NUMBER FIELDS}}
\\
It is very likely that the field defined in Theorem \ref{thm:mproof} 
can be chosen so that $\mbox{h}(K)$ is divisible by $2^{s-1}n.$ 
A sufficient condition for this is that the ideal 
$\sqrt{d}.\OKF+\alpha.\OKF$ is of even order in the 
class group of $KF.$

\begin{tabular}{lr}
 Ivan Chipchakov          & Kalin Kostadinov  \\
 {\it Institute of Mathematics} &  {\it Department of Mathematics} \\
 {\it Bulgarian Academy of Sciences} & {\it Boston University} \\
 {\it Acad G. Bonchev, bl 8}  & {\it 111 Cummington Str.} \\
 {\it 1113 Sofia, Bulgaria}    & {\it Boston, MA02215 USA}   \\
 {\bf chipchak@math.bas.bg} & {\bf kost@math.bu.edu}
 \end{tabular}



\end{document}